\documentclass{amsart}
\usepackage{graphicx}
\usepackage{latexsym}
\usepackage{amsfonts}
\usepackage[all]{xy}
\usepackage{amssymb, mathrsfs, amsfonts, amsmath}
\usepackage{amsbsy}
\usepackage{amsfonts}
\newtheorem*{acknowledgement}{Acknowledgement}
\newtheorem{corollary}{Corollary}
\newtheorem{conjecture}{Conjecture}

\newtheorem{theorem}{Theorem}

\numberwithin{equation}{section}

\begin{document}
\title[4-dimensional compact manifolds]{4-dimensional compact manifolds with nonnegative biorthogonal curvature}
\author{E. Costa$^{1}$ \& E. Ribeiro Jr.$^{2}$}
\address{Universidade Federal da Bahia - UFBA, Departamento de Matem\'{a}tica, Campus de Ondina, Av. Ademar de Barros, 40170-110-Salvador / BA, Brazil.} \email{ezio@ufba.br}
\address{$^{2}$ Universidade Federal do Cear\'a - UFC, Departamento  de Matem\'atica, Campus do Pici, Av. Humberto Monte, Bloco 914,
60455-760-Fortaleza / CE , Brazil.} \email{ernani@mat.ufc.br}
\thanks{$^{2}$ Partially supported by grants from  PJP-FUNCAP/Brazil and CNPq/Brazil}
\date{October 14, 2013}

\begin{abstract}
The goal of this article is to study the pinching problem proposed by S.-T. Yau in 1990 replacing sectional curvature by one weaker condition on biorthogonal curvature. Moreover, we classify 4-dimensional compact oriented Riemannian manifolds with nonnegative biorthogonal curvature. In particular, we obtain a partial answer to Yau Conjecture on pinching theorem for 4-dimensional compact manifolds. 
\end{abstract}

\maketitle
\section{Introduction}
\label{int}

In the last century very much attention has been given to 4-dimensional compact Riemannian manifolds with positive scalar curvature. A classical problem in geometry is to classify  such manifolds in the category of either topology, diffeomorphism, or isometry. This subject have been studied extensively because their connections with general relativity and quantum theory. For comprehensive references on such a theory, we indicate for instance \cite{atiyah}, \cite{besse}, \cite{1}, \cite{DER}, \cite{gursky}, \cite{lebrun}, \cite{MW}, \cite{scorpan} and \cite{7}. Arguably, classifying 4-dimensional compact Riemannian manifolds or understanding their geometry is definitely an important issue.

In 1990 S.-T. Yau collected some important open problems. Here, we call attention to the paragraph where he wrote:
\begin{flushright}
\begin{minipage}[t]{4.37in}
 \emph{``The famous pinching problem says that on a compact simply connected manifold if $K_{min}>\frac{1}{4}K_{max}>0,$ then the manifold is homeomorphic to a sphere. If we replace $K_{max}$ by normalized scalar curvature, can we deduce similar pinching results?"}(cf. \cite{schoen} problem 12, page 369; see also \cite{8}.)
 \end{minipage}
\end{flushright}

In other words, Yau's Conjecture on pinching theorem can be rewritten as follows (cf.  \cite{3}).

\begin{conjecture}[Yau, 1990]
\label{conj1}
Let $(M^n, g)$ be a compact simply connected Riemannian manifold. Denote by $s_{0}$ the normalized scalar curvature of $M^n.$ If $K_{min} > \frac{n-1}{n+2}s_{0},$ then $M^n$ is diffeomorphic to a standard  sphere $\Bbb{S}^n$.
\end{conjecture}

A classical example obtained in \cite{3} shows that $\frac{n-1}{n+2}$ is the best possible pinching for this conjecture (cf. Example 3.1 in \cite{3}). We also notice that if $s$ is the scalar curvature of a  Riemannian manifold $M^n,$ then the normalized scalar curvature of $M^n$ is given by $s_0 = \frac{s}{n(n-1)}.$

From this point on, $M^4$ will denote a compact oriented 4-dimensional manifold and $\mathcal{M}$ the
set of Riemannian metrics $g$ on $M^4$ with  scalar curvature $s$ and  sectional curvature $K.$ Before we state our first theorem, we introduce some definitions. First, let us recall that for each plane $P\subset T_{p}M$ at a point $p\in M^4,$  we define  {\it the biorthogonal (sectional) curvature} of $P$ by the following average of the sectional curvatures
\begin{equation}\label{[1.2]}
K^\perp (P) = \frac{K(P) + K(P^\perp) }{2},
\end{equation}
where $P^\perp$  is the orthogonal plane to  $P.$

The sum of two sectional curvatures on two orthogonal planes  appeared previously in works due to Seaman \cite{6} and Noronha \cite{5}.  It should be remarked that a  compact manifold $M^4$ is Einstein if and only if $K^\perp=K.$ Moreover, $M^4$ is locally conformally flat if and only if $K^{\perp}=s/12.$ We also notice that $\Bbb{S}^{1}\times \Bbb{S}^{3}$ with its canonical metric shows that positive biorthogonal curvature does not imply positive Ricci curvature. Indeed, the positivity of the biorthogonal curvature is an intermediate condition between positive sectional curvature and positive scalar curvature.

Next, we recall that the biorthogonal curvature of a Riemannian manifold $M^4$ is called  {\it weakly $1/4$-pinched} if there exists a positive function $f\in C^{\infty}(M)$ satisfying a suitable pinching condition involving the biorhogonal curvature. In  \cite{6}, Seaman showed that this pinching condition implies nonnegative isotro\-pic curvature. While the first author observed in  \cite{2}  that if $K_1^\perp$ is the minimum of the biorthogonal curvature in each point, then $12K_1^\perp$ is a modified scalar curvature with corresponding modified Yamabe invariant $Y_{1}^{\perp}(M).$ In particular, Costa used the notion of biorthogonal curvature to show a relationship between these invariants and Hopf's Conjecture. We recall that Hopf's Conjecture asks if $\Bbb{S}^2 \times \Bbb{S}^2$ admits a metric with positive sectional curvature.  Costa was able to show Hopf's Conjecture provided  $Y_{1}^{\perp}(\Bbb{S}^2 \times \Bbb{S}^2)\leq 0;$ see \cite{2}. However, Bettiol proved that $Y_{1}^{\perp}(\Bbb{S}^2 \times \Bbb{S}^2)>0,$ which implies that $\Bbb{S}^2 \times \Bbb{S}^2$ admits metrics of positive biorthogonal curvature; for more details see Theorem 1 in \cite{bettiol}. In particular, Bettiol showed that the connected sum $\Bbb{CP}^{2}\# \overline{\Bbb{CP}}^{2}$ admits metrics with positive biorthogonal curvature.

To fix notations we now consider for each point $p\in M^4$ the following functions
\begin{equation}
\label{[1.3]}
K_1^\perp(p) = \textmd{min} \{K^\perp(P); P  \textmd{ is a 2- plane in } T_{p}M \},
\end{equation}

\begin{equation}
\label{[1.4]}
 K_3^\perp(p) = \textmd{max}\{K^\perp (P); P \textmd{ is a 2- plane in } T_{p}M \}
\end{equation}
and
\begin{equation}
\label{[1.5]}
K_2^\perp (p)= \frac{s(p)}{4} - K_1^\perp(p) - K_3^\perp(p).
\end{equation} These functions appeared previously in \cite{2}. In the next section we will collect some properties of these preceding functions.  It is perhaps worth mentioning that the canonical metrics  of the manifolds $\Bbb{S}^4,$ $\Bbb{CP}^2$ and $\Bbb{S}^1\times \Bbb{S}^3$  have $K_1^\perp =s/12,$ $K_1^\perp = s/24$ and $K_1^\perp = s/12,$ respectively.

Our aim is to investigate  the pinching problem on 4-dimensional compact manifolds replacing sectional curvature by biorthogonal curvature conditions. To do this, we start by replacing the assumption $K > s/24$ of the sectional curvature in Conjecture \ref{conj1} by a weaker condition on biorthogonal curvature. With this setting we now announce our first result.

 \begin{theorem}
 \label{thm1}
 Let $(M^4,\, g)$ be a compact oriented Riemannian manifold with positive scalar curvature $s$ satisfying $K_1^\perp \geq s/24.$  Then one of the following assertions holds:
\begin{enumerate}
\item  $(M^4,\,g)$ is diffeomorphic to a connected sum  $\Bbb{S}^4 \sharp(\Bbb{R} \times \Bbb{S}^3)/G_1\sharp ...\sharp(\Bbb{R} \times \Bbb{S}^3)/G_n$, where each $G_i$ is a discrete subgroup of the isometry group of $\Bbb{R} \times \Bbb{S}^3;$
\item  or $(M^4,\, g)$ is isometric to a complex projective space  $\Bbb{CP}^2$ with the Fubini-Study metric.
\end{enumerate}
\end{theorem}

It is worth pointing out that Theorem \ref{thm1} remains true replacing the assumption $K_1^\perp \geq s/24$ by $K_{3}^{\perp}\leq s/6;$ this comment will be clarified in the next section.  As an application of Theorem \ref{thm1} we deduce the following result under finite fundamental group hypothesis.

\begin{corollary}
Let $(M^4, g)$ be a compact oriented Riemannian manifold satisfying  $K_1^\perp \geq s/24.$ We assume that $M^4$ has finite fundamental group. Then we have.
\begin{enumerate}
\item  Either $M^4$ is diffeomhorphic to a sphere $\Bbb{S}^{4}.$
\item Or $(M^4, g)$ is isometric to a complex projective space  $\Bbb{CP}^2$ with the Fubini-Study metric.
\end{enumerate}
\end{corollary}

Proceeding, we now remember that  the space of harmonic 2-forms $H^2(M^4,\Bbb{R})$ can be split as $$H^2(M^4,\Bbb{R}) = H^+(M^4,\Bbb{R})\oplus H^-(M^4,\Bbb{R}),$$ where $H^\pm(M^4,\Bbb{R})$ is the space of positive and negative harmonic 2-forms, respectively. Moreover,  the second Betti number $b_2$ of $M^4$ is $b_2 = b^+ + b^-$,  where $b^\pm$ = dim $H^\pm(M^4,\Bbb{R})$. We recall that $M^4$ is called {\it positive definite} (respectively {\it negative definite}) whether  $b^- = 0$ (respectively $b^+ = 0$). Otherwise, $M^4$ will be called {\it indefinite}. According to Donaldson's and Freedman's works (cf. \cite{donaldson} and \cite{freedman}), if $M^4$ is simply connected and definite, then $M^4$ is homeomorphic to sphere $\Bbb{S}^4$, provided $b_2 = 0$ or  $M^4$  is homeomorphic to a connected sum of $b_2$ complex projective spaces $\Bbb{CP}^2\sharp ... \sharp \Bbb{CP}^2.$

In fact, an elegant argument due to Seaman \cite{7} shows that a compact oriented Riemannian manifold $M^4$ with positive sectional curvature admitting a harmonic 2-form of constant length must be definite. Later, this result  was improved by Noronha in \cite{5}, for more details see Theorems 3.5 and 3.6 in quoted article. More precisely, Noronha proved the following result.

\begin{theorem}[Noronha, \cite{5}]
\label{thmnoronha}
     Let $(M^4,\, g)$  be a compact oriented Riemannian manifold with positive scalar curvature. Then the following assertions hold:
    \begin{enumerate}
   \item
    If $(M^4,\, g)$ admits a non trivial harmonic 2-form of constant length and  $K^\perp > 0,$ then $M^4$ is definite.
   \item
    If  $(M^4, \,g)$ admits a non trivial parallel 2-form  and $K^\perp \geq 0,$ then $M^4$ is biholomhorphic to $\Bbb{CP}^2$ or its universal covering $\widetilde{M}$ is isometric to $M_1^2 \times M_2^2,$ where each $M_i^2$ is diffeomorphic to  sphere $\Bbb{S}^2.$
   \end{enumerate}
\end{theorem}

We now recall that a Riemannian manifold $(M^4,\,g)$ is called {\it geometrically formal} if the wedge product of two harmonic forms is again harmonic; for more details, we refer the reader to \cite{bar}. This concept appeared recently in a work due to Kotschick  \cite{kot}, where he showed that for these class of metrics, harmonic forms have constant length (see also Theorems B and C in \cite{bar}). In particular, Kotschick proved that if $M^4$ is formal and has finite fundamental group, then $M^4$ has second Betti number $b_{2}\leq 2.$ We call attention for the following result due to Kotschick (cf. Corollary 3 in \cite{kot}).

 \begin{theorem}[Kotschick, \cite{kot}]
\label{thmKot}
     Let $M^4$  be a compact simply connected manifold. If $M^4$ is formal and it admits a metric (possibly non formal) with nonnegative scalar curvature, then one of the following assertions occurs:
\begin{enumerate}
\item  $M^4$ is homeomorphic to a sphere $\Bbb{S}^4;$
\item $M^4$ is diffeomorphic to a complex projective space $\Bbb{CP}^2;$
\item or $M^4$ is diffeomorphic to a product of two spheres $\Bbb{S}^2 \times \Bbb{S}^2.$
\end{enumerate}
  \end{theorem}

More recently, B\"ar in \cite{bar} was able to prove the following classification under nonnegative sectional curvature assumption.

\begin{theorem}[B\"ar, \cite{bar}]
\label{thmbar}
     Let $(M^4,\, g)$  be a compact oriented geometrically formal  Riemannian manifold.
    \begin{enumerate}
\item  If $M^4$ is simply connected and $(M^4,\,g)$ has sectional curvature $K\geq 0.$ Then:
\begin{enumerate}
\item $M^4$ is homeomorphic to a sphere $\Bbb{S}^4;$
\item $M^4$ is diffeomorphic  to a complex projective space $\Bbb{CP}^2;$
\item or $(M^4,\, g) = M_1^2 \times M_2^2$, where each $M_i^2$ is diffeomorphic to  a sphere $\Bbb{S}^2$ and has nonnegative sectional curvature.
\end{enumerate}
\item If $(M^4,\,g)$ has positive sectional curvature.  Then:
\begin{enumerate}
\item  Either $M^4$ is homeomorphic to a sphere $\Bbb{S}^4;$
\item or $M^4$ is diffeomorphic  to a complex projective space $\Bbb{CP}^2.$
\end{enumerate}
\end{enumerate}
\end{theorem}

As was previously mentioned, we are interested in classfying 4-dimensional manifolds under biorthogonal curvature hypotheses. To this end,  also note that from Micallef and Moore work \cite{4} nonnegative sectional curvature implies  $K_{3}^\perp \leq s/4,$ and nonnegative isotropic curvature implies $K_{3}^\perp \leq s/4.$ For more details see the next section. Based in these observations and  inspired by Seaman, Noronha, Kotschick and B\"ar ideas developed in \cite{7}, \cite{5}, \cite{kot} and \cite{bar}, respectively,  we now  announce our second theorem.

\begin{theorem}
\label{thm2}
Let $(M^4,\, g)$  be a compact oriented Riemannian manifold with positive scalar curvature. Then the following assertions occur:
    \begin{enumerate}
   \item
    If $(M^4,\, g)$ admits a non trivial harmonic 2-form of constant length and  $K_3^\perp < \frac{s}{4},$ then $M^4$ is definite.
   \item
    If  $(M^4,\, g)$ admits a non trivial parallel 2-form  and $K_3^\perp \leq \frac{s}{4},$ then $M^4$ is biolomorphic to a complex projective space $\Bbb{CP}^2$ with the Fubini-Study metric or its universal covering $\widetilde{M}$ is isometric to $M_1^2 \times M_2^2,$ where each $M_i^2$ is homeomorphic to a sphere $\Bbb{S}^2.$
   \end{enumerate}
  \end{theorem}

We point out that  Theorem \ref{thm2} can be seen as an improvement to Theorems \ref{thmnoronha} and \ref{thmbar}. In particular, we obtain the following characterization under an integral condition involving the biorthogonal curvature.

  \begin{corollary}
\label{cor1thm2}
Let $(M^4,\, g)$  be a compact oriented simply connected Riemannian manifold with positive scalar curvature $s$ and satisfying $$\int_{M}\big(s - 4K_{3}^\perp\big)dV_g \geq 0.$$
 If all harmonic forms of $(M^4, g)$ have constant length,  then one of the following assertions occurs:
\begin{enumerate}
\item  $M^4$ is homeomorphic to $\Bbb{S}^4;$
\item $M^4$ is diffeomorphic to $\Bbb{CP}^2;$
\item  or $M^4$ is isometric to $M_1^2 \times M_2^2,$  where each $M_i^2$ is diffeomorphic to a sphere $\Bbb{S}^2.$
\end{enumerate}
     \end{corollary}

In order to state the next result, we adopt the following notation. For an oriented manifold $M^4,$ we consider $\Lambda^2$ be the bundle of 2-forms $\alpha\in M^4$ and let $\ast: \Lambda^2\to \Lambda^2$ be the Hodge star operator. Thus, there is a invariant decomposition $\Lambda^2=\Lambda^{+}\oplus\Lambda^{-},$ where $\Lambda^{\pm}=\{\alpha\in\Lambda^{2};\,\ast\alpha=\pm\alpha\},$ depending only on the orientation and the conformal class of the metric. Therefore, the Weyl curvature tensor $W$ is an endomorphism of  $\Lambda^2$ such that $W = W^+\oplus W^-,$ where    $W^\pm : \Lambda^\pm
\longrightarrow \Lambda^\pm$ are called of the self-dual and anti-self-dual parts of $W.$ Half conformally flat metrics are also known as self-dual or anti-self-dual if $W^{-}$ or $W^{+}=0,$ respectively. These metrics are, in a certain sense, analogous to anti-self-dual connections in Yang-Mills theory.

The formal divergence $\delta$ for any  $(0,4)$-tensor $T$ is given by  $$\delta T(X_1,X_2,X_3) =
-trace_{g}\{(Y,Z)\mapsto\nabla_{Y}T(Z,X_1,X_2,X_3)\},$$ where $g$ is
the metric of $M^4$. We say that the Weyl tensor of $M^4$ is harmonic when  $\delta W=0.$

One fundamental inequality in Riemannian geometry is {\it Kato's inequality}. Namely, let $s\in\Gamma(E),$ where $E\to M$ is a vector bundle over $M,$ then  $|\nabla|s||\leq |\nabla s|$ away from the zero locus of $s.$
In a famous article LeBrun and Gursky proved a {\it refined  Kato's inequality}. More precisely, they showed that if $W^{+}$ is harmonic, then away from the zero locus of $W^{+}$ we have
\begin{equation}
\label{katoinequality}
|d|W^{+}||\leq \sqrt{\frac{3}{5}}|\nabla W^{+}|.
\end{equation} Moreover, (\ref{katoinequality}) holds in the distributional sense on $M^4,$ for more details see Lemma 2.1 in \cite{gursky}.

On the basis of these observations and inspired on ideas developed in \cite{yang} we use improved Kato's inequality jointly with a classical theorem due to Hitchin \cite{besse} to prove the following result.

\begin{theorem}
\label{thm3}
Let $(M^4,\, g)$  be a compact oriented Riemannian manifold  with harmonic Weyl tensor and positive scalar curvature. We assume that $g$ is analytic  and  $$K_{1}^{\perp}\geq \frac{s^{2}}{8(3s+5\lambda_{1})},$$ where $\lambda_{1}$ stands for  the first eigenvalue of the Laplacian with respect to $g.$ Then one of the following assertions holds:
\begin{enumerate}
\item  
$M^4$ is diffeomorphic to a connected sum  $\Bbb{S}^4 \sharp(\Bbb{R} \times \Bbb{S}^3)/G_1\sharp ...\sharp(\Bbb{R} \times \Bbb{S}^3)/G_n$, where each $G_i$ is a discrete subgroup of the isometry group of $\Bbb{R} \times \Bbb{S}^3.$ In this case, $g$ is locally conformally flat;
\item  or $M^4$ is isometric to a complex projective space  $\Bbb{CP}^2$ with the Fubini-Study metric.
\end{enumerate}
\end{theorem}

As an immediate consequence, we obtain the following corollary.

\begin{corollary}
\label{cor1thm3}
Let $(M^4,\, g)$  be a compact oriented Riemannian manifold  with harmonic Weyl tensor and metric $g$ analytic. We assume that $Ric\geq \rho >0$ and  $K_{1}^{\perp}\geq \frac{3s^{2}}{8(9s^{2}+20\rho)},$  then we have. 
\begin{enumerate}
\item Either $M^4$ is isometric to $\Bbb{S}^4$ with its canonical metric.
\item Or $M^4$ is isometric to $\Bbb{CP}^{2}$ with Fubini-Study metric.
\end{enumerate}
\end{corollary}

We already know that a compact manifold $M^4$ is Einstein if and only if $K^\perp=K.$  Moreover, from Theorem 5.26 in \cite{besse} Einstein metrics are analytic. It should be emphasized that there are regularity results which could be used to show that harmonic self-dual Weyl tensor implies that 
the metric is analytic choosing appropriate coordinates (e.g. harmonic one), for more details see \cite{DeTurck}. Finally, we deduce the following corollary which was first obtained by Yang in \cite{yang}.

\begin{corollary}
Let $(M^4,\, g)$  be a compact oriented Einstein manifold  satisfying $Ric= \rho >0.$ Suppose $$K \geq \frac{2\rho^2}{12\rho+5\lambda_1}.$$
 Then either $M^4$ is isometric to $\Bbb{S}^4$ with its canonical metric or $M^4$ is isometric to $\Bbb{CP}^{2}$ with Fubini-Study metric.
\end{corollary}

\section{Preliminaries}

Throughout this section we collect a couple of formulae that will be useful in the proofs of our results. As was previously commented the Weyl tensor $W$ is an endomorphism of
$\Lambda^2 $ such that $W = W^+\oplus W^-,$ where    $W^\pm : \Lambda^\pm\longrightarrow \Lambda^\pm$ are called the self-dual and anti-self-dual parts of $W,$ respectively. Furthermore, if $\mathcal{R}$ denotes the curvature of $M^4$ we get a matrix

\begin{equation}
\mathcal{R}=
\left(
  \begin{array}{c|c}
    \\
W^{+} +\frac{s}{12}Id & B \\ [0.4cm]\hline\\

    B^{*} & W^{-}+\frac{s}{12}Id  \\[0.4cm]
  \end{array}
\right),
\end{equation}
where $B:\Lambda^{-}\to \Lambda^{+}$ stands for the Ricci traceless operator of $M^4$ given by $B=Ric-\frac{s}{4}g.$ For more details in this subject, we recommend the famous ``Besse's book"  \cite{besse}.

We now consider  $w_1^\pm \leq w_2^\pm \leq w_3^\pm$ be the eigenvalues of the tensors $W^\pm,$ respectively. In \cite{2}, the first author proved some formulae involving the biorthogonal curvatures and the eigenvalues of $W^\pm$ that will be important in the proofs of our results. Here we present their proofs for the sake of completeness.

First, we consider a point $p\in M^4$ and $X,Y\in T_{p}M$ orthonormal. Therefore, the unitary 2-form $\alpha=X\wedge Y$ can be uniquely written as $\alpha=\alpha^{+}+\alpha^{-},$ where $\alpha^{\pm}\in \Lambda^{\pm}$ with $|\alpha^{+}|^{2}=\frac{1}{2}$ and $|\alpha^{-}|^{2}=\frac{1}{2}.$ Moreover, under these notations the sectional curvature $K(\alpha)$ is given by
\begin{equation}
\label{k}
K(\alpha)=\frac{s}{12}+\langle\alpha^{+},W^{+}(\alpha^{+})\rangle+\langle\alpha^{-},W^{-}(\alpha^{-})\rangle+2\langle\alpha^{+},B\alpha^{-}\rangle.
\end{equation}
In particular, we have
\begin{equation}
\label{kperp}
K(\alpha^{\perp})=\frac{s}{12}+\langle\alpha^{+},W^{+}(\alpha^{+})\rangle+\langle\alpha^{-},W^{-}(\alpha^{-})\rangle-2\langle\alpha^{+},B\alpha^{-}\rangle,
\end{equation} where $\alpha^{\perp}=\alpha^{+}-\alpha^{-}.$  Combining (\ref{k}) with (\ref{kperp}) we arrive at
\begin{equation}
\frac{K(\alpha)+K(\alpha^{\perp})}{2}=\frac{s}{12}+\langle\alpha^{+},W^{+}(\alpha^{+})\rangle+\langle\alpha^{-},W^{-}(\alpha^{-})\rangle.
\end{equation}
Hence, we can use (\ref{[1.3]}) to obtain
\begin{eqnarray*}
\label{eq898}
K_{1}^{\perp}=\frac{s}{12}+min\left\{\langle\alpha^{+},W^{+}(\alpha^{+})\rangle; \,|\alpha^{+}|^{2}=\frac{1}{2}\right\}+min\left\{\langle\alpha^{-},W^{-}(\alpha^{-})\rangle; \,|\alpha^{-}|^{2}=\frac{1}{2}\right\}.
\end{eqnarray*}
However, from Proposition 2.1 of \cite{noronha2}  there exists an orthonormal basis of $\Lambda^{2}$ given by $$\{X_{1}\wedge Y_{1},\,X_{2}\wedge Y_{2},X_{3}\wedge Y_{3}\},$$ where $X_{i},\,Y_{i}\in T_{p}M$ for all $i=1,2,3.$ In particular, we invoke Proposition 2.5 also in \cite{noronha2} to get
\begin{equation}
\label{[1.6]}
K_1^\perp - \frac{s}{12} = \frac{w_1^+ + w_1^-}{2}.
\end{equation}
Arguing in the same way we obtain
\begin{equation}\label{[1.8]}
K_3^\perp - \frac{s}{12} = \frac{w_3^+ + w_3^-}{2}.
\end{equation}
Finally, from (\ref{[1.5]}) we have
\begin{equation}\label{[1.7]}
K_2^\perp - \frac{s}{12} = \frac{w_2^+ + w_2^-}{2}.
\end{equation}

From results due to Micallef and Moore $M^4$ has nonnegative isotropic curvature if and only if $w_3^\pm \leq s/6;$ for more details see \cite{4}. For that reason, nonnegative sectional curvature implies  $K_{3}^\perp \leq s/4,$ and nonnegative isotropic curvature implies $K_{3}^\perp \leq s/4.$  Moreover, we notice that $K \geq s/24$ implies that $K_1^\perp \geq s/24,$  as well as  $K_1^\perp \geq s/24$ implies that $K_3^\perp\leq s/6.$

\section{Proof of the results}

\subsection{The Proof of Theorem \ref{thm1}}

\begin{proof}
Let $(M^4,\, g)$ be a compact oriented Riemannian manifold with positive scalar curvature $s.$ Since $K_1^\perp \geq \frac{s}{24}$ implies $K_3^\perp \leq \frac{s}{6}$ it is enough to assume   that  $K_3^\perp \leq s/6.$ Now, from (\ref{[1.8]}) we arrive at $$w_3^+ + w_3^- =2K_3^\perp - \frac{s}{6}.$$ From this, we may use that $w_1^\pm \leq w_2^\pm \leq w_3^\pm$ and $w_1^\pm + w_2^\pm+ w_3^\pm=0$ jointly with our assumption to conclude that   $w_3^{+} \leq s/6.$ Similarly, we conclude that $w_3^{-}\leq s/6,$ which implies that $M^4$ has nonnegative isotropic curvature.  Assume that $M^4$ admits a metric with positive isotropic curvature, then $M^4$  is diffeomorphic to a connected sum  $\Bbb{S}^4 \sharp(\Bbb{R} \times \Bbb{S}^3)/G_1\sharp ...\sharp(\Bbb{R} \times \Bbb{S}^3)/G_n$, where each $G_i$ is a discrete subgroup of the isometry group of $\Bbb{R} \times \Bbb{S}^3.$

On the other hand, we assume that $M^4$ does not admit a metric with positive isotropic curvature. Hence, if $M^4$ is irreducible, we can apply Theorem 1.1 of \cite{8} to deduce that $(M^4,\,g)$ is either locally symmetric or is K\"ahler. We now suppose that $M^4$ is irreducible and locally symmetric, which implies that $(M^4,\,g)$ is an Einstein manifold. Therefore, we may use Theorem 4.4 of \cite{MW} to conclude that $(M^4,\,g)$ is isometric to a complex projective space $\Bbb{CP}^2.$ In the K\"ahler case it is known that $w_3^+ = \frac{s}{6}.$ For this, we invoke (\ref{[1.8]}) to obtain $$w_3^-  \leq -\frac{s}{6} + 2K_3^\perp - \frac{s}{6} \leq 0,$$  which implies that $w_3^-  = 0$ in $M^4.$ From this it follows that  $W^- = 0.$ Now, we apply Theorem 1.1 in \cite{DER} to conclude that $(M^4,\, g)$ is locally symmetric and then $(M^4,\, g)$ is isometric to a complex projective space $\Bbb{CP}^2$.

Finally, we consider  $(M^4,\, g)$ locally reducible. Since  $K_3^\perp \leq \frac{s}{6},$ it is not difficult to check that this case  can not occur;  for more details see  Theorem 3.1 in \cite{MW}. So, we conclude the proof of the theorem.
\end{proof}

\subsection{The Proof of Theorem \ref{thm2}}

 The first part of the proof will follow from Noronha (cf. Theorem 3.6 in \cite{5}). First of all we consider $(M^4,\, g)$ be a compact oriented  Riemannian manifold with positive scalar curvature $s.$ Moreover, let  $ \alpha^+\in \Lambda^+$  (respectively $\alpha^-\in \Lambda^-$)  be a positive (respectively negative) non-degenerate differentiable 2-form.  From this, we have two Weitzenb\"och formulae given by
\begin{equation}
\label{[1.9]}
\langle \Delta\alpha^\pm, \alpha^\pm \rangle = \frac{1}{2} \Delta\mid  \alpha^\pm \mid ^2 + \mid   \nabla \alpha^\pm   \mid^2 + \langle(\frac{s}{3} - 2W^\pm)\alpha^\pm, \alpha^\pm \rangle.
\end{equation}

We now denote $w_3^\pm$ be the largest eigenvalues of $W^\pm$, respectively. Under these conditions we have $$\langle W^{\pm}(\alpha^{\pm}),\alpha^{\pm}\rangle\leq w_{3}^{\pm}\langle\alpha^{\pm},\alpha^{\pm}\rangle.$$ It then follows from (\ref{[1.9]}) that
\begin{equation}
\label{1.10}
\langle \Delta\alpha^\pm,  \alpha^\pm \rangle  \geq  \frac{1}{2} \Delta \mid \alpha^\pm  \mid^2 +  \mid \nabla \alpha^\pm   \mid^2 + \left(\frac{s}{3} - 2w_3^\pm\right)   \mid\alpha^\pm   \mid^2.
\end{equation}

Now we are ready to prove Theorem \ref{thm2}.

\begin{proof} First, we assume that $M^4$ has a non-degenerate harmonic 2-form $\alpha$ with constant length. Suppose that $M^4$ is not definite. This means that  $b^\pm > 0,$ which gives the following possibilities:

 \begin{enumerate}
\item $\alpha$ is a negative 2-form.
\item $\alpha$ is a positive 2-form.
\item $\alpha = \alpha^+ \ + \ \alpha^-$, where $\alpha^\pm$ are non-degenerate positive and negative 2-form, respectively.
\end{enumerate}

We suppose that occurs the first case (the second case has similar argue). Thus, we may use (\ref{1.10}) for $\alpha = \alpha^-$
 to deduce $$0 \geq \   \mid \nabla \alpha  \mid^2 + \left(\frac{s}{3} - 2w_3^-\right)  \mid\alpha  \mid^2.$$ From this it follows that $w_3^- \geq \frac{s}{6}$ in $M^4.$ Next, from (\ref{[1.8]}) we have $$w_3^- + w_3^+ = 2K_3^\perp - \frac{s}{6} <\frac{s}{3}$$ and then $w_3^+ + w_3^- <\frac{s}{3},$ which implies $w_{3}^{+}<\frac{s}{6}.$

 On the other hand, insofar as $b^+ > 0$, there exists a harmonic non-degenerate positive 2-form $\gamma.$ Furthermore,  (\ref{1.10}) with respect to  $\gamma$  ensures
$$0\geq \frac{1}{2}\Delta |\gamma|^{2}+|\nabla \gamma|^{2}+\left(\frac{s}{3}-2w_{3}^{+}\right)|\gamma |^{2}.$$ We integrate the last expression and we use Stokes theorem to obtain
$$0 \geq  \int_{M}\mid   \nabla \gamma \mid^2dV_g + \int_{M}\left(\frac{s}{3} - 2w_3^{+}\right) \mid\gamma \mid^2dV_g > 0,$$ which is a contradiction. This proves the first possibility.

We now treat the third case. To this end, we use that $\alpha = \alpha^+ \ + \ \alpha^-$ jointly with (\ref{1.10}) to infer
\begin{eqnarray*}
0&\geq& |\nabla \alpha^{+}|^{2}+|\nabla\alpha^{-}|^{2}+\left(\frac{s}{3}-2w_{3}^{+}\right)|\alpha^{+}|^{2}\\&&+\left(\frac{s}{3}-2w_{3}^{-}\right)|\alpha^{-}|^{2}.
\end{eqnarray*}
After a straightforward computation we get
\begin{eqnarray}
\label{[1.11]}
0 &\geq& |\nabla\alpha^{+}|^{2}+|\nabla\alpha^{-}|^{2}+\left(s-4K_{3}^{\perp}\right)|\alpha^{+}|^{2}\nonumber\\&&+\left(2w_{3}^{-}-\frac{s}{3}\right) \left(|\alpha^{+}|^{2}-|\alpha^{-}|^{2}\right).
\end{eqnarray}

Next, if there exists a point $p\in M^{4}$ such that $|\alpha^{+}|^{2}=|\alpha^{-}|^{2}$ we use (\ref{[1.11]}) to deduce $$0\geq|\nabla\alpha^{+}|^{2}+|\nabla\alpha^{-}|^{2}+\left(s-4K_{3}^{\perp}\right)|\alpha^{+}|>0,$$ which is again a contradiction. Therefore, without loss of generality we can assume that $|\alpha^{+}|^{2}>|\alpha^{-}|^{2}.$ From this it follows that $w_{3}^{-}\leq s/6$ in $M^4.$

On the other hand, on integrating (\ref{1.10}) for $\alpha^{-}$ we obtain
\begin{equation}
\label{[1.12]}
0\geq \int_{M}|\nabla\alpha^{-}|^{2}dV_{g}+\int_{M}\left(\frac{s}{3}-2w_{3}^{-}\right)|\alpha^{-}|^{2}dV_{g}\geq 0,
\end{equation} which implies that $\nabla\alpha^{-}=0,$ and then $|\alpha^{-}|$ is constant. By using once more (\ref{[1.12]}) we conclude that $w_{3}^{-}=s/6$ and  whereas $w_{3}^{+}+w_{3}^{-}< s/3,$ we infer $w_{3}^{+}<s/6.$ Finally, we take the integral in (\ref{1.10}) to get $$0\geq \int_{M}|\nabla\alpha^{+}|^{2}dV_{g}+\int_{M}\left(\frac{s}{3}-2w_{3}^{+}\right)|\alpha^{+}|^{2}dV_{g}>0,$$ which is once more a contradiction. So, we have proved the first assertion of the theorem.

Proceeding we suppose that $(M^4,\,g)$ admits a non trivial parallel 2-form. Under this condition it is well-known that $M^4$ is K\"ahler, in particular, $w_{3}^{+}=\frac{s}{6}.$ Since $K_3^\perp \leq \frac{s}{4},$ we conclude that $w_{3}^{-}\leq s/6$ and then from Micallef-Moore work $M^4$ has nonnegative isotropic curvature (cf. \cite{4}). Next, if $M^4$ is locally irreducible we use Theorem 1.2 due to Seshadri \cite{8} to conclude that $M^4$ is biholomorphic to $\Bbb{CP}^2$ or isometric to a compact Hermitian symmetric space. In the last case, $M^4$ is Einstein and then it is isometric to $\Bbb{CP}^2.$ To finish, it suffices to argue as in the proof of Theorem \ref{thm1} to deduce that $M^4$ can not be locally reducible. \end{proof}

\subsection{The Proof of Corollary \ref{cor1thm2}}
\begin{proof} We assume that $M^4$ is simply connected and that their harmonic forms have constant length. We now assume the unpublished theorem (at present) due to Kotschick (cf. Theorem \ref{thmKot}) to deduce that $M^4$ is either homeomorphic to a sphere $\Bbb{S}^4,$ diffeomorphic to a complex projective space $\Bbb{CP}^2$ or diffeomorphic to  product of two spheres $\Bbb{S}^2 \times \Bbb{S}^2.$ In the last case, we have the second Betti number $b_2^\pm = 1.$ Therefore, we consider two harmonic forms with constant length $\alpha^\pm \in H^\pm(M, \Bbb{R}),$ without loss of generality, we may assume that the length is equal to one.  Whence, on integrating (\ref{[1.9]}) we obtain
 $$0 \geq \int_{M}  \mid \nabla \alpha^+  \mid^2dV_g + \int_{M}  \mid \nabla \alpha^-  \mid^2dV_g + \int_{M} \big( s - 4K_3^\perp\big)dV_g \geq 0.$$
From this it follows that $\alpha^\pm$ are parallel and by using once more (\ref{[1.9]})  for $\alpha^\pm$  we infer $w^\pm = s/6$ and  $K_3^\perp = s/4.$ Finally, we invoke Theorem \ref{thm2} to conclude the proof of the corollary.
\end{proof}

\subsection{The Proof of Theorem \ref{thm3}}
\begin{proof}
Since $\delta W=0$ we have the following Weitzenb\"ock formulae (cf. 16.73 in \cite{besse})
\begin{equation}
\label{Weit}
\frac{1}{2}\Delta |W^{\pm}|^{2}+|\nabla W^{\pm}|^{2}+\frac{s}{2}|W^{\pm}|^{2}-18 \det W^{\pm} =0.
\end{equation}
Moreover, by use of Lagrange multipliers we infer
\begin{equation}
\label{lagrange}
\det W^{\pm}\leq \frac{\sqrt{6}}{18}|W^{\pm}|^3.
\end{equation}
However, our hypothesis implies that $|W^{\pm}|^{2}$ are analytic. So far, the set $$\Sigma=\{p\in M;\,|W^{+}|(p)=0\,\,\hbox{or}\,\,|W^{-}|(p)=0\}$$ is finite provided $W^{\pm}\not\equiv 0.$ Suppose by contradiction that $(M^4,\,g)$ is not half conformally flat. For this, there is a constant $t>0$ such that
\begin{equation*}
\int_{M}\big(|W^{+}|-t|W^{-}|\big)dV_{g}=0.
\end{equation*}
Choosing $W^{-}$ in  (\ref{Weit}) and multiplying by $t^{2}$ and adding the result to (\ref{Weit}) with respect to $W^{+},$ we deduce
\begin{eqnarray}
\label{eq3thm3}
0&\geq& \big(|\nabla W^{+}|^{2}+t^{2}|\nabla W^{-}|^{2}\big)+\frac{s}{2}\big(|W^{+}|^{2}+t^{2}|W^{-}|^{2}\big)\nonumber\\&&-18\big(\det W^{+}+t^{2}\det W^{-}\big).
\end{eqnarray}
Applying  refined  Kato's inequality (\ref{katoinequality}) jointly with (\ref{lagrange}) in the previous inequality we get
\begin{eqnarray}
\label{12309}
0&\geq& \frac{5}{3}\int_{M}\big(|d|W^{+}||^{2}+t^{2}|d|W^{-}||^{2}\big)dV_{g}+\int_{M}\frac{s}{2}\big(|W^{+}|^{2}+t^{2}|W^{-}|^{2}\big)dV_{g}\nonumber\\&&-\sqrt{6}\int_{M}\big(|W^{+}|^{3}+t^{2}|W^{-}|^{3}\big)dV_{g}.
\end{eqnarray}
On the other hand, we notice that
\begin{eqnarray*}
\big(|d|W^{+}||^{2}+t^{2}|d|W^{-}||^{2}\big)&=&\frac{1}{2}\big(|d(|W^{+}|-t|W^{-}|)|^{2}+|d(|W^{+}|+t|W^{-}|)|^{2}\big)\\&\geq& \frac{1}{2}|d(|W^{+}|-t|W^{-}|)|^{2}.
\end{eqnarray*}
Moreover, from Poincar\'e Inequality we have
\begin{equation*}
\frac{1}{2}\int_{M}|d(|W^{+}|-t|W^{-}|)|^{2} dV_{g}\geq \frac{\lambda_{1}}{2}\int_{M}(|W^{+}|-t|W^{-}|)^{2} dV_{g},
\end{equation*}
so that, from the two preceding inequalities we obtain
\begin{equation}
\label{eq5thm3}
\int_{M}\big(|d|W^{+}||^{2}+t^{2}|d|W^{-}||^{2}\big)dV_{g}\geq \frac{\lambda_{1}}{2}\int_{M}(|W^{+}|-t|W^{-}|)^{2} dV_{g}.
\end{equation}
Therefore, comparing (\ref{eq5thm3}) with (\ref{12309}) we obtain
\begin{eqnarray*}
0&\geq& \frac{5}{6}\lambda_{1}\int_{M}\big(|W^{+}|^{2}-2t|W^{+}||W^{-}|+t^{2}|W^{-}|^{2}\big)dV_{g}
\nonumber\\&&+\int_{M}\frac{s}{2}\big(|W^{+}|^{2}+t^{2}|W^{-}|^{2}\big)dV_{g}-\sqrt{6}\int_{M}\big(|W^{+}|^{3}+t^{2}|W^{-}|^{3}\big)dV_{g},
\end{eqnarray*}
which can be written as
\begin{eqnarray}
\label{7676}
0&\geq&\int_{M}\Big\{|W^{-}|^{2}\Big(\frac{5}{6}\lambda_{1}+\frac{s}{2}-\sqrt{6}|W^{-}|\Big)t^{2} -\Big(\frac{5}{3}\lambda_{1}|W^{+}||W^{-}|\Big) t \nonumber\\&&+|W^{+}|^{2}\Big(\frac{5}{6}\lambda_{1}+\frac{s}{2}-\sqrt{6}|W^{+}|\Big)\Big\}dV_{g}.
\end{eqnarray}
For simplicity, we can write the integrand of (\ref{7676}) as
\begin{equation}
\label{ghgh}
\mathcal{P}(t)= |W^{-}|^{2}\left(a-\sqrt{6}|W^{-}|\right)t^{2}-\frac{5}{3}\lambda_{1}|W^{+}||W^{-}| t+ |W^{+}|^{2}\left(a-\sqrt{6}|W^{+}|\right),
\end{equation}
where $a=\frac{5}{6}\lambda_{1}+\frac{s}{2}.$ We notice that (\ref{ghgh}) is a quadratic function of $t$ and its discriminant $\Delta$ is given by
\begin{equation}
\label{Delta}
\Delta = \frac{25}{9}\lambda_{1}^{2}|W^{+}|^{2}|W^{-}|^{2}-4|W^{+}|^{2}|W^{-}|^{2}\big(a-\sqrt{6}|W^{+}|\big)\big(a-\sqrt{6}|W^{-}|\big).
\end{equation}
On the other hand, we recall that $|W^{\pm}|\leq 6(w_{1}^{\pm})^{2}$  and then we use (\ref{[1.6]}) to deduce  $$|W^{+}|+|W^{-}|\leq \sqrt{6}\left(\frac{s}{6}-2K_{1}^{\perp}\right).$$ A straightforward computation shows that our assumption on biorthogonal curvature implies $$\sqrt{6}\left(\frac{s}{6}-2K_{1}^{\perp}\right)\leq \frac{4a^{2}-\frac{25}{9}\lambda_{1}^{2}}{4\sqrt{6}a}.$$ Whence,
\begin{equation}
\label{pop}
|W^{+}|+|W^{-}|\leq \frac{4a^{2}-\frac{25}{9}\lambda_{1}^{2}}{4\sqrt{6}a}.
\end{equation}
Therefore, we combine (\ref{pop}) with (\ref{Delta})  to conclude that $\Delta$ is less than or equal to zero. Hence, we use once more  (\ref{7676})  to deduce $|W^{+}||W^{-}|=0$ in $M^4$. But, since  $\Sigma$ is finite we arrive at a contradiction. 

Therefore, we conclude that $W^{+}=0$ or $W^{-}=0.$ Finally, we define the following sets $$A=\big\{p\in M^{4};\,Ric(p)\neq\frac{s(p)}{4}g\big\}$$ and $$B=\big\{p\in M^{4};\,|W^{+}|(p)=|W^{-}|(p)\big\},$$ where $(Ric-\frac{s}{4}g)$ stands for the traceless Ricci tensor of $(M^{4},\,g).$ If $A$ is empty we conclude that $M^4$ is an Einstein manifold. In this case, we invoke Hitchin's theorem \cite{besse} to conclude that $M^4$ is either isometric to $\Bbb{S}^4$ with its canonical metric or isometric to $\Bbb{CP}^{2}$ with Fubini-Study metric. Otherwise, if $A$ is not empty, then there exists  a point $p\in M^4$ and an open set $U$ such that $p\in U\subset A.$ So, we use Corollary 1 of \cite{derd2} to conclude that $U\subset A\subset B.$ So far, since the function $f=|W^{+}|^{2}-|W^{-}|^{2}$ is analytic we conclude that $f$ is identically zero. For this, $|W^{+}|^{2}=|W^{-}|^{2}$ and then $M^4$ is locally conformally flat. This implies that $M^4$ has positive isotropic curvature and then we use once more Chen-Tang-Zhu theorem  \cite{1} to conclude that $M^4$ is diffeomorphic to a connected sum  $\Bbb{S}^4 \sharp(\Bbb{R} \times \Bbb{S}^3)/G_1\sharp ...\sharp(\Bbb{R} \times \Bbb{S}^3)/G_n$, where each $G_i$ is a discrete subgroup of the isometry group of $\Bbb{R} \times \Bbb{S}^3.$ This finishes the proof of the theorem.
\end{proof}

\subsection{Proof of the Corollary \ref{cor1thm3}}

\begin{proof}
Since $Ric\geq \rho >0$ implies $\lambda_{1}\geq \frac{4\rho}{3}$ we combine Theorem \ref{thm3}  with Tani's theorem in \cite{T} and then gives the promised result.
\end{proof}

\begin{acknowledgement}
The authors want to thank the referees for their careful reading and helpful suggestions. The second author would like to thank the Institute of Mathematics - UFBA where part of this work was carried out. He wish to express his gratitude for the excellent support during his stay. \nonumber
\end{acknowledgement}

\end{document}